\documentclass[12pt]{amsart}

\usepackage[mathscr]{eucal}
\usepackage{amscd}
\usepackage{amsfonts}
\usepackage{amsmath, amsthm, amssymb}
\usepackage{latexsym}
\usepackage[dvips]{graphics}

\theoremstyle{plain} 
\newtheorem{thm}{Theorem}

\theoremstyle{definition}
\newtheorem{defn}[thm]{Definition}
\newtheorem{ex}[thm]{Example}

\theoremstyle{remark}

\numberwithin{equation}{section}

\usepackage[all,poly,necula]{xy}
\newcommand{\M}{\mathbb{M}}
\newcommand{\pos}{\mathbb{Z}_{>0}}

\newsavebox\locboxinminipage
\newlength\locboxinminipagel
\newcommand{\boxinminipage}[1]
{%
 \sbox\locboxinminipage{#1}%
 \settowidth\locboxinminipagel{\usebox{\locboxinminipage}}%
 \begin{minipage}{\locboxinminipagel}\usebox{\locboxinminipage}\end{minipage}%
}
\DeclareMathOperator{\Hom}{Hom}

\title{INTRODUCTION TO ALGEBRAS OF PARTIAL TRIANGULATIONS}      

\author{Laurent Demonet} 

\address{
\begin{flushleft}
        \hspace{0.3cm}  Graduate School of Mathematics \\
         \hspace{0.3cm}  Nagoya University\\
         \hspace{0.3cm}  Furocho, Chikusaku \\
         \hspace{0.3cm}  Nagoya 464-8602 JAPAN\\
\end{flushleft}
}
\email{Laurent.Demonet@normalesup.org}

\thanks{The paper is in a final form and no version of it will
be submitted for publication elsewhere.}

\begin{document}



\begin{abstract}

The aim of this note is to give a gentle introduction to algebras of partial triangulations of marked surfaces, following the structure of a talk given during the 49th symposium on ring theory and representation theory, held in Osaka. This class of algebras, which always have finite rank, contains classical Jacobian algebras of triangulations of marked surfaces and Brauer graph algebras. We discuss representation theoretical properties and derived equivalences. All results are proven in \cite{Laurent-Demonet:2}, under slightly milder hypotheses.




\end{abstract}

\maketitle 

\section{The algebra of a partial triangulation}

Let $k$ be a unital ring and $\Sigma$ be a connected compact oriented surface with or without boundary. We fix a non-empty finite set $\M \subset \Sigma$ of marked points (some of them may be on the boundary $\partial \Sigma$). For each $M \in \M$, we fix an invertible coefficient $\lambda_M \in k^\times$ and a multiplicity $m_M \in \pos$. For simplicity, we suppose here that if $\Sigma$ is a sphere then $\# \M \ge 5$ and if $\Sigma$ is a disc then $\# \M \ge 3$. 

\begin{defn}
 An \emph{arc} on $(\Sigma, \M)$ is a continuous map $u: [0, 1] \to \Sigma$ satisfying:
 \begin{itemize}
  \item The restriction of $u$ to $(0,1)$ is an embedding into $\Sigma \setminus \M$;
  \item Extremities $0$ and $1$ are mapped to $\M$.
 \end{itemize}
 We consider arcs up to homotopy relative to their endpoints in $\Sigma \setminus \M$. Moreover, we exclude arcs that are homotopic to a marked point or to a boundary component, that is the closure of a connected component of $\partial \Sigma \setminus \M$.
\end{defn}

In this note, for simplicity, we exclude arcs that are loops enclosing a unique marked point $M \in \M$ with $m_M \leq 2$. We do not make this restriction in \cite{Laurent-Demonet:2}.

\begin{defn}
 We say that two arcs $u$ and $v$ are \emph{compatible} if, up to homotopy, they are non-crossing. Then, a \emph{partial triangulation of $(\Sigma, \M)$} is a set $\sigma$ of arcs of $(\Sigma, \M)$ that are pairwise compatible. 
  If $\sigma$ is a maximal partial triangulation and each connected component of $\partial \Sigma$ contains at least a marked point, $\sigma$ is called a \emph{triangulation}.
\end{defn}

In order to define the algebra $\Delta_\sigma$ of a partial triangulation $\sigma$, we first need to construct a quiver $Q_\sigma$. 

\begin{defn}
 The quiver $Q_\sigma$ has set of vertices $\sigma$, and arrows are winding in $\Sigma$ between successive arcs around marked points counter-clockwisely. We call \emph{bouncing path} of $Q_\sigma$ any path of length $2$ consisting of two arrows that are not successive around the same endpoint.
\end{defn}

 We give two examples to illustrate this definition:

\begin{ex} \label{expt} In the left example $\Sigma$ is a disc with four marked points $A$, $B$, $C$ and $D$ ($A$ and $B$ are on the boundary). In the right one, $\Sigma$ is a torus with two marked points $M$ and $N$. The partial triangulations are depicted with thick lines.
 $$
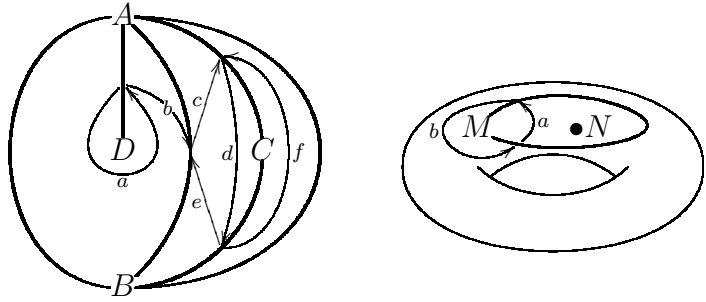
\xymatrix@L=.05cm@R=1.5cm@C=1.5cm@M=0.01cm{
   & A \ar@{-}@`{c+(-20,0),p+(-20,0)}[dd]|(.25){}^{{}}="AB" \ar@{-}@`{c+(-20,0),p+(-20,0)}[dd]<-.3pt> \ar@{-}@`{c+(-20,0),p+(-20,0)}[dd]<.3pt> \ar@{-}[d]|(.25){}^{{}}="AD" \ar@{-}[d]<-.3pt> \ar@{-}[d]<.3pt> \ar@{-}@`{c+(+35,0),p+(+35,0)}[dd] \ar@{-}@`{c+(+35,0),p+(+35,0)}[dd]<-.3pt> \ar@{-}@`{c+(+35,0),p+(+35,0)}[dd]<.3pt>\\
   & D & C \ar@{-}@`{c+(0,10),p+(10,0)}[ul]|(.25){}^{{}}="CA" \ar@{-}@`{c+(0,10),p+(10,0)}[ul]<-.3pt> \ar@{-}@`{c+(0,10),p+(10,0)}[ul]<.3pt> \\
   & B \ar@{-}@`{c+(10,0),p+(0,-10)}[ur]|(.25){}^{{}}="BC" \ar@{-}@`{c+(10,0),p+(0,-10)}[ur]<-.3pt> \ar@{-}@`{c+(10,0),p+(0,-10)}[ur]<.3pt> \ar@{-}@`{c+(12,10), p+(12,-10)}[uu]|(.25){}^{{}}="BA" \ar@{-}@`{c+(12,10), p+(12,-10)}[uu]<-.3pt> \ar@{-}@`{c+(12,10), p+(12,-10)}[uu]<.3pt>
   \ar@/^/|b"AD";"BA"
   \ar^c"BA";"CA"
   \ar@`{c+(-8,-8),p+(0,-16),p+(8,-8)}_a"AD";"AD"
   \ar^e"BC";"BA"
   \ar@/^/_d"CA";"BC"
   \ar@`{c+(12,0),p+(12,0)}_f"BC";"CA"
 } \quad \quad \quad \xymatrix@L=.05cm@!=0cm@R=.5cm@C=.5cm@M=0.00cm{
   & \\
   & \\
   & \\
   & & M \ar@{-}@`{c+(5,-4),c+(35,3),p+(5,8)}[]_{}|(.12)*+<.1cm>{{}}="su"|(.88){{}}="tu" \ar@{-}@`{c+(5,-4),c+(35,3),p+(5,8)}[]<.3pt> \ar@{-}@`{c+(5,-4),c+(35,3),p+(5,8)}[]<-.3pt> & & & \bullet N\\
   \ar@{-}@`{c+(0,15),p+(0,15)}[rrrrrrrr] \ar@{-}@`{c+(0,-15),p+(0,-15)}[rrrrrrrr] & & \ar@{-}@`{c+(5,-5),p+(-5,-5)}[rrrr]|(.1){{}}="P"|(.9){{}}="Q" &  & & & & & \\
   & \\
   & 
   \ar@{-}@`{c+(5,4),p+(-5,4)}"P";"Q"
   \ar@/_/"su";"tu"_a
   \ar@`{c+(-5,0),c+(-10,-3),p+(-5,-4)}"tu";"su"_b
 }$$

 Bouncing paths of the first example are $ab$, $cd$, $de$ and $ec$. Bouncing paths of the second example are $a^2$ and $b^2$.
\end{ex}

To define the algebra $\Delta_\sigma$, we need two more combinatorial concepts:
\begin{defn}
 For each edge $u \in \sigma$ and endpoint $M$ of $u$, we denote by $\omega_{u, M}$ the path of $Q_\sigma$ going from $u$ to $u$ winding once around $M$ if $M \notin \partial \Sigma$ and $\omega_{u, M} = 0$ if $M \in \partial \Sigma$. 
\end{defn}

\begin{defn}
 A \emph{small triangle} of $\sigma$ is a triple $(u, v, w)$ of arcs of $\sigma$ together with three marked points $M$, $N$ and $P$ such that 
  \begin{itemize}
   \item $u$ is incident to $M$ and $N$;
   \item $v$ is incident to $N$ and $P$;
   \item $w$ is incident to $P$ and $M$;
   \item The union of $u$, $v$ and $w$ encloses clockwisely (in the order $u$, $v$, $w$) a disc without marked point inside.
  \end{itemize}
\end{defn}

We now define the \emph{algebra of the partial triangulation $\sigma$} by $\Delta_\sigma := k Q_\sigma / I_\sigma$ where $I_\sigma$ is the ideal generated by the following relations:
\begin{enumerate}
 \item For each $u \in \sigma$ that joins $M \in \M$ to $N \in \M$, we require $\lambda_M \omega_{u, M}^{m_M} = \lambda_N \omega_{u, N}^{m_N}$.
 \item For each $u \in \sigma$ with endpoint $M$, we require $\omega_{u, M}^{m_M} q = 0$ for any arrow $q$ of $Q_\sigma$.
 \item Suppose that $\sigma$ contains a small triangle $M, N, P$ as in the following picture:
 $$
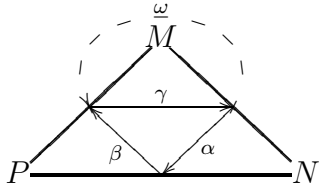
\xymatrix@L=.05cm@R=1.5cm@C=1.5cm@M=0.01cm{
  & M \ar@{-}[dl]|{}="MP" \ar@{-}[dl]<-.3pt> \ar@{-}[dl]<.3pt> \ar@{-}[dr]|{}="MN" \ar@{-}[dr]<-.3pt> \ar@{-}[dr]<.3pt> &  \\
  P \ar@{-}[rr]|{}="NP" \ar@{-}[rr]<-.3pt> \ar@{-}[rr]<.3pt> & & N
  \ar"MP";"MN"^\gamma
  \ar"MN";"NP"^\alpha
  \ar"NP";"MP"^\beta
  \ar@{-->}@`{c+(8,16),p+(-8,16)}"MN";"MP"_\omega
  }$$
  where $\alpha$, $\beta$ and $\gamma$ are arrows, $\omega$ is a the only possible simple path winding around $M$ if $M \notin \partial \Sigma$ and $\omega = 0$ if $M \in \partial \Sigma$. We require the relation $\alpha \beta = \lambda_M (\omega \gamma)^{m_M - 1} \omega$.
 \item For any bouncing path $\alpha \beta$ that does not appear in case 3, we require $\alpha \beta = 0$.
\end{enumerate}

\begin{ex}
 Relations for the left partial triangulation of Example \ref{expt} are:
 \begin{align*}
  a^{m_D} &= (df)^{m_C} = (fd)^{m_C} = 0 & (1) \\
  ab &= 0 & (4) \\
  cd &= de = 0 & (3) \\
  ec &= \lambda_C (fd)^{m_C - 1} f & (3)
 \end{align*}
 and other relations are redundant.

 Relations for the right partial triangulation are: 
 \begin{align*}
  (ab)^{m_M} &= (ba)^{m_M} & (1) \\
  a^2 &= b^2 = 0. & (4)
 \end{align*}
\end{ex}

It turns out that algebras of partial triangulations are particularly well behaved. A first result about them is that this definition is compatible with the naive notion of a sub-partial triangulation:
\begin{thm} \label{thmsubtri}
 Let $\tau \subset \sigma$. Then we have 
 $$\Delta_\tau \cong e_\tau \Delta_\sigma e_\tau$$
 where $e_\tau$ is the sum of the primitive idempotents of $\Delta_\sigma$ corresponding to the arcs of $\tau$. 
\end{thm}
Notice that we have naturally $k Q_\tau \subset k Q_\sigma$. However, relations as defined in this note do not go through this inclusion. We have to take a more complicated variant of these relations, giving an isomorphic algebra, to obtain Theorem \ref{thmsubtri}.

\section{Brauer graph algebras and Jacobian algebras of surfaces}

We explain here that the class of algebras of partial triangulations contains two important classes of algebras.

\begin{thm} \label{thbr}
 If $\sigma$ contains no small triangle, neither arc incident to $\partial \Sigma$, then $\Delta_\sigma$ is the \emph{Brauer graph algebra} of $\sigma$ considered as a \emph{ribbon graph}.
 Moreover, any Brauer graph algebra is the algebra of a partial triangulation of a surface without boundary.
\end{thm}
We will not recall what a Brauer graph algebra is. For more details, see for example \cite{Laurent-Demonet:4} or \cite{Laurent-Demonet:5}. However, this definition is very close to the definition of the algebra of a partial triangulation and Theorem \ref{thbr} is mostly straightforward.

\begin{thm}
 If all $m_M$ are invertible in $k$ and $\sigma$ is a triangulation, then $\Delta_\sigma$ is the Jacobian algebra of a quiver with potential $(Q_\sigma, W_\sigma)$. To define $W_\sigma$, consider the set $E$ of small triangles $T$ of $\sigma$ up to rotation and for each of them denote by $\alpha_T$, $\beta_T$ and $\gamma_T$ the three arrows as in the figure defining relation (3) earlier. Then, for each $M \in \M$ take arbitrarily an arc $u_M$ incident to $M$. Then
 $$W_\sigma := \sum_{T \in E} \alpha_T \beta_T \gamma_T - \sum_{M \in \M} \frac{\lambda_M}{m_M}  \omega_{u, M}^{m_M}$$
 (as usual for potentials, terms are only well defined up to cyclic permutations).
\end{thm}
 Notice that if $m_M = 1$ for all $M \in \M$, we recover the usual Jacobian algebra of a surface as defined in \cite{Laurent-Demonet:3}. Recall that the Jacobian algebra of $(Q_\sigma, W_\sigma)$ is the quotient of the completed path algebra $\widehat{k Q_\sigma}$ by the cyclic derivatives of $W_\sigma$. So we get here an improvement as $\Delta_\sigma$ is directly defined from $k Q_\sigma$ without completion.

\section{Algebraic properties of $\Delta_\sigma$}

 We have the following result about $\Delta_\sigma$:
\begin{thm} \label{thmdim}
 The $k$-algebra $\Delta_\sigma$ is a free $k$-module of rank
 $$\sum_{M \in \M \cap \partial \Sigma} \frac{d_M(d_M-1)}{2} + \sum_{M \in \M \cap (\Sigma \setminus \partial \Sigma)} m_M d_M^2 + f$$
 where, for $M \in \M$, $d_M$ is the degree of $M$ in the graph $\sigma$ (without counting boundary components), and $f$ is the number of arcs in $\sigma$ with both endpoints on boundaries. 
\end{thm}

\begin{ex}
 The algebra of the left partial triangulation of Example \ref{expt} has rank $5 + 4 m_C + m_D$, and the right one has rank $4 m_M$.
\end{ex}

More precisely, there is a $k$-basis of $\Delta_\sigma$ consisting of all strict and non-idempotent prefixes of all $\omega_{u, M}^{m_M}$, together with primitive idempotents and elements $\lambda_M \omega_{u, M}^{m_M} = \lambda_N \omega_{u, N}^{m_N}$.

The following property generalizes a known result for Brauer graph algebras and Jacobian algebras of surfaces without boundary:
\begin{thm}
 If $\sigma$ has no arc incident to the boundary, then $\Delta_\sigma$ is a symmetric $k$-algebra (\emph{i.e.} $\Hom_k(\Delta_\sigma, k) \cong \Delta_\sigma$ as $\Delta_\sigma$-bimodules). 
\end{thm}

\section{Representation type of $\Delta_\sigma$}

The next theorem permits to expect that the classification of $\Delta_\sigma$-modules is possible:
\begin{thm}
 If $k$ is an algebraically closed field, then $\Delta_\sigma$ is of tame representation type.
\end{thm}
The proof of this result relies on a deformation theorem by Crawley-Boevey \cite{Laurent-Demonet:1}. Indeed, the relations defining $\Delta_\sigma$ can be deformed to the relations of a Brauer graph algebra in a suited manner. Moreover, Brauer graph algebras are of tame representation type. Notice that unfortunately, these techniques do not permit to deduce directly the classification of $\Delta_\sigma$-modules even though modules over Brauer graph algebras are known.

\section{Flip of partial triangulations and derived equivalences}

Finally, we give a flip leading to derived equivalences. For an arc $u$ in $\sigma$ such that arcs marked by $+$ in the following diagrams are also in $\sigma$ (in particular, they are not in $\partial\Sigma$), we define $\mu_u(\sigma)$ by replacing $u$ by $u^*$ defined in the following way:
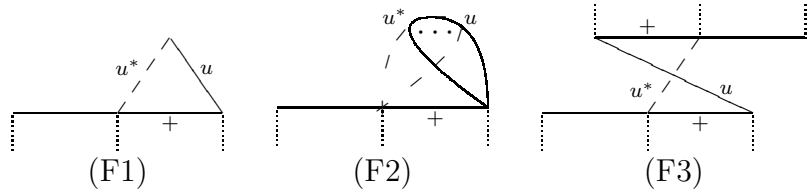
\begin{center} \begin{tabular}{ccccc}
                 \boxinminipage{$$\xymatrix@L=.05cm@!=0cm@R=.5cm@C=.7cm@M=0.01cm{
                  & & & \\
		    & & &  \\
			  & & & & \\
		    \ar@{..}[d] & & \ar@{--}[uur]^{ u^*} \ar@{-}[ll] \ar@{..}[d]  & & \ar@{-}[uul]_{ u} \ar@{-}[ll]^{+} \ar@{..}[d]  \\ & & & & 
                 }$$} &  &
		\boxinminipage{$$\xymatrix@L=.05cm@!=0cm@R=.5cm@C=.7cm@M=0.01cm{
                     & & &  \cdots \\
			   & & & \\
		   \ar@{..}[d] & & \ar@{-}@`{p+(14,10),p+(8,14),p+(0,10)}@{--}[]_(.7){ u^*} \ar@{..}[d] \ar@{-}[ll] & & \ar@{-}@`{p+(0,10),p+(-8,14),p+(-14,10)}_(.3){ u} \ar@{-}[ll]^{+} \ar@{..}[d]  \\ & & & & 
                 }$$} &  & 
                \boxinminipage{$$\xymatrix@L=.05cm@!=0cm@R=.5cm@C=.7cm@M=0.01cm{
                  & & & & & \\
		    & \ar@{..}[u] \ar@{-}[rr]^{+}  & & \ar@{..}[u] \ar@{-}[rr] & & \ar@{..}[u]  \\ & & & & 
			  & \\
		    \ar@{..}[d] & & \ar@{--}[uur]^(.2){u^*} \ar@{-}[ll] \ar@{..}[d]  & & \ar@{-}[uulll]_(.2){ u} \ar@{-}[ll]^{+} \ar@{..}[d]  \\ & & & & 
                 }$$} \\
	    (F1) & &  (F2) & & (F3)
                \end{tabular}  
 \end{center}
 We also define coefficients $\mu_u(\lambda)$. For a marked point $M \in \M$, $\mu_u(\lambda)_M = \lambda_M$ except in the following cases:
 \begin{itemize}
  \item In case (F1), if $M$ is the topmost vertex of the figure then $\mu_u(\lambda)_M = -\lambda_M$.
  \item In case (F1), if $M$ is the rightmost vertex of the figure then $\mu_u(\lambda)_M = (-1)^{m_M}\lambda_M$.
  \item In case (F2), if $M$ is the unique marked point enclosed by $u$ then $\mu_u(\lambda)_M = -\lambda_M$.
 \end{itemize}

 Then we get the following result:
 \begin{thm}
  There is a derived equivalence between $\Delta_\sigma$ and $\Delta_{\mu_u(\sigma)}^{\mu_u(\lambda)}$ where the second algebra is computed with respect to the coefficients $\mu_u(\lambda)$.
 \end{thm}

\begin{ex}
 We consider the two following partial triangulations of a disc with three marked points, none of them are in $\partial\Sigma$:
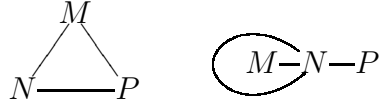
 \begin{center}
  \boxinminipage{$$\xymatrix@L=.05cm@!=0cm@R=1cm@C=.7cm@M=0.01cm{
                  & M \ar@{-}[dl] \ar@{-}[dr] & \\
                  N \ar@{-}[rr] & & P
                 }$$} \quad \quad \quad 
  \boxinminipage{$$\xymatrix@L=.05cm@!=0cm@R=.5cm@C=.7cm@M=0.01cm{
                   \\
                  M \ar@{-}[r] & N \ar@{-}[r] \ar@{-}@`{p+(-7,7), p+(-20,0), p+(-7,-7)}[] & P  \\
                 }$$}  
 \end{center}
 They are related by a flip so the following algebras, obtained for $\lambda_M = \lambda_N = \lambda_P = 1$ and $m_M = m_N = m_P = m$ are derived equivalent:
 $$\frac{k \left(\boxinminipage{\xymatrix@L=.05cm@!=0cm@R=2cm@C=1.4cm{
        & \bullet \ar@/_/[dr]_x \ar@/_/[dl]_y \\
        \bullet \ar@/_/[ur]_x \ar@/_/[rr]_y & & \bullet \ar@/_/[ul]_y \ar@/_/[ll]_x}}\right)}{(x^2-(yx)^{m-1}y, y^2)}
  \quad \text{and} \quad
   \frac{k \left(\boxinminipage{\xymatrix@L=.05cm@!=0cm@R=1cm@C=1.2cm{
            \\
         & \bullet \ar@/_/[r]_{\beta_1} \ar@`{p+(-5,5), p+(-15,0), p+(-5,-5)}[]_\alpha & \bullet \ar@/_/[r]_{\gamma_1} \ar@/_/[l]_{\beta_2} & \bullet \ar@/_/[l]_{\gamma_2} \ar@`{p+(5,-5), p+(15,0), p+(5,5)}[]_\delta & \\  & }}
        \right)}{\left(\begin{array}{c} \beta_2 \alpha - (\gamma_1 \gamma_2 \beta_2 \beta_1)^{m-1} \gamma_1 \gamma_2 \beta_2, \\ \alpha \beta_1 - (\beta_1 \gamma_1 \gamma_2 \beta_2)^{m-1} \beta_1 \gamma_1 \gamma_2, \\ \beta_1 \beta_2 - \alpha^{m-1}, \gamma_1 \delta, \delta \gamma_2, \gamma_2 \gamma_1,\\ \delta^m - (\gamma_2 \beta_2 \beta_1 \gamma_1)^m \end{array}\right)}.
$$

\end{ex}


 
\ifx\undefined\bysame 
\newcommand{\bysame}{\leavevmode\hbox to3em{\hrulefill}\,} 
\fi 


\end{document}